\documentclass[a4paper,10pt,reqno]{amsart}
\usepackage{graphicx}
\usepackage{geometry}
\usepackage{float}
\usepackage{verbatim}
\usepackage{setspace}
\usepackage{siunitx}
\usepackage{paralist}
\usepackage{enumitem}
\usepackage{xcolor}
\usepackage{amsfonts,amsmath,amssymb,bbm,amsthm,amsbsy}
\usepackage[utf8]{inputenc}
\usepackage{lipsum}
\usepackage[english]{babel}
\usepackage{booktabs}
\usepackage{array}
\usepackage{subfig}
\usepackage{caption}
\usepackage{mathtools}
\usepackage{stmaryrd}
\usepackage{hyperref}
\usepackage{stackengine,scalerel}

\newenvironment{frcseries}{\fontfamily{frc}\selectfont}{}

\makeatletter
\newtheorem*{rep@theorem}{\rep@title}
\newcommand{\newreptheorem}[2]{
\newenvironment{rep#1}[1]{
 \def\rep@title{#2 \ref{##1}}
 \begin{rep@theorem}}
 {\end{rep@theorem}}}
\makeatother

\newtheorem{thm}{Theorem}[section]
\newtheorem{lemma}[thm]{Lemma}
\newtheorem{prop}[thm]{Proposition}
\newtheorem{corr}[thm]{Corollary}

\newtheorem{innercustomthm}{Theorem}
\newenvironment{customthm}[1]
  {\renewcommand\theinnercustomthm{#1}\innercustomthm}
  {\endinnercustomthm}

\newtheorem{innercustomprop}{Proposition}
\newenvironment{customprop}[1]
  {\renewcommand\theinnercustomprop{#1}\innercustomprop}
  {\endinnercustomprop}

\newtheorem{innercustomcorr}{Corollary}
\newenvironment{customcorr}[1]
  {\renewcommand\theinnercustomcorr{#1}\innercustomcorr}
  {\endinnercustomcorr}

\newtheorem*{thm*}{Theorem}
\newtheorem*{lemma*}{Lemma}
\newtheorem*{prop*}{Proposition}
\newtheorem*{corr*}{Corollary}
\newtheorem*{claim*}{Claim}

\theoremstyle{remark}
\newtheorem{rmk}[thm]{Remark}

\newtheorem*{rmk*}{Remark}
\newtheorem*{conj*}{Conjecture}
\newtheorem*{quest*}{Question}

\theoremstyle{definition}
\newtheorem{defn}[thm]{Definition}
\newtheorem{exmp}[thm]{Example}

\newtheorem*{defn*}{Definition}
\newtheorem*{exmp*}{Example}

\newreptheorem{theorem}{Theorem}
\newreptheorem{corollary}{Corollary}
\newreptheorem{proposition}{Proposition}

\usepackage{fancyhdr}
\usepackage{tikz}
\usepackage{tikz-cd}

\usepackage{IEEEtrantools}

\newenvironment{equ*}[1]{\begin{IEEEeqnarray*}{#1}}{\end{IEEEeqnarray*}}

\newcommand{\F}{\mathbb{F}}

\newcommand{\Z}{\mathbb{Z}}

\newcommand{\Nc}{\mathcal{N}}

\newcommand{\C}{\mathcal{C}}

\DeclareFontFamily{U}{mathx}{}
\DeclareFontShape{U}{mathx}{m}{n}{<-> mathx10}{}
\DeclareSymbolFont{mathx}{U}{mathx}{m}{n}
\DeclareMathAccent{\widecheck}{0}{mathx}{"71}

\newcommand{\inj}{\hookrightarrow}
\newcommand{\sur}{\twoheadrightarrow}

\newcommand{\RG}{R\llbracket G\rrbracket}

\DeclareMathOperator{\coplus}{\,\widehat{\oplus}}
\newcommand{\oprod}{\ensurestackMath{\stackinset{c}{}{c}{-.33pt}%
  {\boldsymbol{\cdot}}{\displaystyle\prod}}}
\DeclareMathOperator{\bigcoplus}{\,\widehat{\bigoplus}}
\DeclareMathOperator{\crotimes}{\,\widehat{\otimes}_\mathit{R}\,}

\DeclareMathOperator{\Fd}{\textup{Fin.}\,\textup{dim}}

\DeclareMathOperator{\Ext}{Ext}

\DeclareMathOperator{\Tor}{Tor}

\DeclareMathOperator{\Ind}{Ind}

\DeclareMathOperator{\Res}{Res}

\newcommand{\Top}{\mathbf{Top}}

\newcommand{\Pro}{\mathbf{Pro}}

\newcommand{\PModR}{\mathbf{PMod}(R)}
\newcommand{\coshR}{\mathbf{CoSh}(R)}
\newcommand{\coshZ}{\mathbf{CoSh}(\widehat{\Z})}
\newcommand{\shR}{\mathbf{Sh}(R)}
\newcommand{\DModR}{\mathbf{DMod}(R)}

\newcommand{\M}{\mathcal{M}}

\title{Profinite Direct Sums with Applications to Profinite Groups of Type $\Phi_R$}
\author{Jiacheng Tang}
\thanks{
\\Email: jiacheng.tang@postgrad.manchester.ac.uk; Address: Alan Turing Building, University of Manchester, Manchester, M13 9PY, United Kingdom.}

\calclayout
\begin{document}
\maketitle

\begin{abstract}
We show that the ``profinite direct sum" is a good notion of infinite direct sums for profinite modules, having properties similar to those of direct sums of abstract modules. For example, the profinite direct sum of projective modules is projective, and there is a Mackey's Formula for profinite modules described using these sums. As an application, we prove that the class of profinite groups of type $\Phi_R$ is closed under subgroups.
\end{abstract}

\section{Introduction}
\label{sec1}

Let $R$ be a profinite ring. The category of profinite $R$-modules is abelian with well-behaved limits, but categorical coproducts are not exact in general. This means that, in practice, coproducts of profinite modules are difficult to work with. Instead, \cite{melfree} and \cite{ribesgraph} propose a different kind of direct sum, which we call the \emph{profinite direct sum} (following \cite{gareth}), that is much better-behaved. For example, \cite[Chapter 9]{ribesgraph} uses this profinite direct sum to establish a Mayer-Vietoris sequence for profinite groups acting on profinite trees. Since its invention in \cite{melfree}, the profinite direct sum has been treated mostly algebraically, until a categorical formalisation was found independently in \cite{gareth} and \cite{boggi}\footnote{We will mostly follow the approach of \cite{gareth}.}. This connection will prove useful in the development of the theory.

\subsection*{Aim and Main Results}

The aim of this paper is to show that the profinite direct sum has excellent module-theoretic and homological properties. This makes profinite direct sums significantly more useful than categorical coproducts of profinite modules in practice. The following two results give an indication of why that is the case. We point out the strong resemblance to abstract modules over an abstract ring.

Given a family of profinite $R$-modules $\{M_x\}_{x\in X}$ indexed by a profinite space $X$ (with some compatibility conditions), we denote their profinite direct sum by $\bigcoplus_{x\in X}M_x$ or simply $\bigcoplus_XM_x$, which is also a profinite $R$-module (Definition \ref{fakesh}). The reader should think of this as being analogous to the usual direct sum $\bigoplus_XM_x$ of abstract modules $M_x$ indexed by a set $X$.

\begin{customcorr}{\ref{impcorr}\ref{impcorr2}} The profinite direct sum $\bigcoplus_XM_x$ has projective dimension $\leq d$ if and only if for each $x\in X$, $M_x$ has projective dimension $\leq d$.
\end{customcorr}

\begin{customthm}{\ref{mackey}}[Profinite Mackey's Formula]
Let $G$ be a profinite group and $H,K$ closed subgroups (which are not necessarily open). If $M$ is a profinite right $R\llbracket H\rrbracket$-module, then there is an isomorphism of profinite right $R\llbracket K\rrbracket$-modules $$\Res^G_K\Ind^G_HM=\bigcoplus_{g\in H\backslash G/K}\Ind^K_{K\cap g^{-1}H g}\Res^{g^{-1}H g}_{K\cap g^{-1}H g}Mg.$$
\end{customthm}

Although these results seem purely algebraic, \cite{gareth} proves an equivalence between profinite direct sums (algebraic viewpoint) and ``spaces of global cosections of cosheaves" (categorical viewpoint), which we state formally in Proposition \ref{propff}. The ability to switch between the two perspectives will prove to be extremely valuable, as proofs are often easier from one perspective than the other.

The results stated above should be general enough that they may be used as tools in future works on the cohomology theory or representation theory of profinite groups (such as works similar to \cite{boggi} or \cite{johnmodular}). As an illustration of how they could be useful, we highlight a concrete application of profinite direct sums/Mackey's Formula to profinite groups \emph{of type $\Phi_R$}. To avoid clogging the introduction with technical background, we will simply define type $\Phi_R$ here and leave the motivations to the start of Section \ref{sec3.1}.

Let $G$ be a profinite group and $R$ a profinite ring. We say that $G$ is \emph{of type $\Phi_R$} (Definition \ref{phidef}) if a profinite $\RG$-module has finite projective dimension whenever its restriction to every finite subgroup of $G$ has finite projective dimension. For some of the proofs to work, the ring $R$ needs to satisfy some conditions, which we will state in Section \ref{sec3}, but examples include $R=\F_p$ and $R=\Z_p$.

\begin{customcorr}{\ref{corrtype}}
Let $G$ be a profinite group. Suppose one of the following holds:
\begin{enumerate}[label=(\alph*)]
\item $G$ has finite $vcd_R$.
\item $G$ acts on a profinite tree $T$ with finite (vertex and edge) stabilisers.
\end{enumerate}
Then $G$ is of type $\Phi_R$.
\end{customcorr}

\begin{customthm}{\ref{thmeg}}
There exist profinite groups of type $\Phi_{\F_p}$ which are not of finite $vcd_{\F_p}$.
\end{customthm}

\begin{customprop}{\ref{phisub}}
If a profinite group $G$ is of type $\Phi_R$ then so is every closed subgroup $H$ of $G$.
\end{customprop}

The proofs of the results above will use facts about profinite direct sums that we develop along the way.

\subsection*{Outline of Paper}

We will assume basic knowledge of profinite groups (refer to \cite{profinite}), profinite graphs (refer to \cite{ribesgraph}), category theory (refer to \cite{maccat}), sheaf theory (refer to \cite{macsheaves}) and homological algebra (refer to \cite{weibel}). Profinite graphs are used only in the last subsection (Section \ref{sec3.3}) to construct examples, so readers who are unfamiliar with the topic should not be deterred. Homological algebra of profinite modules (or more generally of pseudocompact modules) can be found in \cite{brumer} and \cite{compactmodules}.

In Section \ref{sec2.1}, we define sheaves and cosheaves over profinite rings (categorical viewpoint) and explain their connections to profinite direct sums (algebraic viewpoint). The ability to switch between the two viewpoints will allow us to easily establish categorical properties of cosheaves (Section \ref{sec2.2}). In Appendix \ref{appendix}, we also construct a tensor product for cosheaves. Section \ref{sec2.3} contains our main results, including Mackey's Formula (Theorem \ref{mackey}).

The results we prove in Section \ref{sec2} are then applied in Section \ref{sec3} to understand profinite groups of type $\Phi_R$. Section \ref{sec3.2} is devoted to exploring the properties of profinite groups of type $\Phi_R$, whilst the reader can find examples of such groups in Section \ref{sec3.3}, which include all profinite groups of finite $vcd_R$ (with some conditions on $R$).

Remark: Whenever we write ``=" in this paper, we mean isomorphic, usually canonically isomorphic (or equivalent in the case of categories).

Convention: All rings are associative with a 1 but are not necessarily commutative (except in Section \ref{sec3}). When unspecified, all of our one-sided profinite modules are right modules.

\subsection*{Acknowledgements}
The author would like to thank his supervisor Peter Symonds for his constant guidance, the following people (in alphabetical order) for helpful discussions on the subject matter: Matthew Antrobus, Calum Hughes, Gregory Kendall, Adrian Miranda, Gareth Wilkes, and Julian Wykowski, as well as the following people for reading drafts of this paper and giving useful feedback: Matthew Antrobus, Gregory Kendall. The author is immensely grateful to the referees for giving multiple suggestions that improved the writing and for pointing out several references that the author had failed to include.

This work was supported by the EPSRC Doctoral Training Partnership [grant number \linebreak EP/W524347/1].

\subsection*{Remarks on v3 and Statement of AI Use}
The arXiv version v3 is to fix some mistakes found by ChatGPT 5.6 Sol. Note that these mistakes are present in the published version. The main mistake is the (previous) proof of Proposition \ref{corrected} and there are also minor issues with Corollary \ref{corrtype}\ref{corra}.

\section{(Co)sheaves over Profinite Rings}
\label{sec2}

In this section, we shall define sheaves and cosheaves over profinite rings and study their properties. The first two subsections are more categorical, so readers looking for algebraic applications might want to move swiftly on to Section \ref{sec2.3}, where we discuss properties of the profinite direct sum $\coplus\,$. The only crucial definition is that of the profinite direct sum, which can be found in Definition \ref{fakesh}.

\subsection{Definitions and Equivalence}
\label{sec2.1}

In this subsection, we will give the categorical definition of cosheaves (Definition \ref{defsh}) and the algebraic viewpoint (Definition \ref{fakesh}) and explain how they are equivalent.

Let $\Top$ and $\Pro$ denote the categories of topological spaces and profinite spaces respectively. Given a profinite ring $R$, let $\DModR$ and $\PModR$ denote the categories of discrete (left) and profinite (right) $R$-modules respectively, which are dual by Pontryagin Duality. Given a topological space $X$, let $O_c(X)$ denote the poset category of clopen subspaces of $X$.

\begin{defn}[\cite{gareth}]\label{defsh}
Let $R$ be a profinite ring. A \emph{presheaf of discrete (topological) $R$-modules} is a pair $(\M,X)$, where $X$ is a profinite space and $\M$ is a functor $O_c(X)^{\mathrm{op}}\to\DModR$. A \emph{morphism of presheaves} $(\Phi,f)\colon(\M,X)\to(\Nc,Y)$ consists of a continuous map $f\colon Y\to X$ and a natural transformation $\Phi\colon\M\circ p_1\to \Nc\circ f^*$ between the functors
\begin{eqnarray*}
\M\circ p_1\colon O_c(X)^{\mathrm{op}}\times O_c(Y)^{\mathrm{op}}\to\DModR,&&(U,V)\mapsto \M(U),\\
\Nc\circ f^*\colon O_c(X)^{\mathrm{op}}\times O_c(Y)^{\mathrm{op}}\to\DModR,&&(U,V)\mapsto \Nc( f^{-1}(U)\cap V),
\end{eqnarray*}
where $p_1$ is projection onto the first component.

A presheaf $(\M,X)$ is a \emph{sheaf} if for every clopen subspace $U\subseteq X$ and every finite clopen cover $U=U_1\cup\cdots\cup U_n$, the following canonical diagram is an equaliser: $$\M(U)\to\prod_{i}\M(U_i)\rightrightarrows\prod_{i,j}\M(U_i\cap U_j).$$ By \cite[Proposition 1.4]{gareth}, this is equivalent to saying that for every clopen subspace $U\subseteq X$ and every decomposition $U=U_1\sqcup\cdots\sqcup U_n$ into finitely many clopen subspaces, the natural map $\M(U)\to \prod^n_{i=1}\M(U_i)$ is an isomorphism.

Let us denote the category of sheaves of discrete $R$-modules by $\mathbf{Sh}(R)$.

Dually, a \emph{precosheaf of profinite $R$-modules} is a pair $(\M,X)$, where $X$ is a profinite space and $\M$ is a functor $O_c(X)\to\PModR$. A \emph{morphism of precosheaves} $(\Psi,f)\colon(\M,X)\to(\Nc,Y)$ consists of a continuous map $f\colon X\to Y$ and a natural transformation $\Psi\colon\M\circ f^*\to\Nc\circ p_2$ between the functors
\begin{eqnarray*}
\M\circ f^*\colon O_c(X)\times O_c(Y)\to\PModR,&&(U,V)\mapsto \M(U\cap f^{-1}(V)),\\
\Nc\circ p_2\colon O_c(X)\times O_c(Y)\to\PModR,&&(U,V)\mapsto \Nc(V),
\end{eqnarray*}
where $p_2$ is projection onto the second component.

A precosheaf $(\M,X)$ is a \emph{cosheaf} if for every clopen subspace $U\subseteq X$ and every decomposition $U=U_1\sqcup\cdots\sqcup U_n$ into finitely many clopen subspaces, the natural map $\bigoplus^n_{i=1}\M(U_i)\to\M(U)$ is an isomorphism.

Let us denote the category of cosheaves of profinite $R$-modules by $\mathbf{CoSh}(R)$.

(Sometimes, we will abuse notation and write $\M$ for a (co)sheaf instead of $(\M,X)$ when $X$ is understood.)
\end{defn}

The above definition seems at first glance to be quite complicated, but later we will give a simpler, equivalent definition (see Proposition \ref{defsame} below).

By \cite[Proposition 1.5]{gareth}, the cosheaf category $\mathbf{CoSh}(R)$ is dual to the sheaf category $\shR$ via Pontryagin Duality. Although sheaves are perhaps more common than cosheaves in the literature, many of our results later will be from the perspective of cosheaves, since that was the viewpoint when cosheaves of profinite modules were first considered (see \cite{melfree}).

As pointed out in \cite[Section 1]{gareth}, the main differences between $\shR$ and usual sheaf categories are that the space $X$ can vary, and that the sheaf condition talks about finite clopen covers (rather than arbitrary open covers). Note that, by compactness of profinite spaces, we can change ``finite clopen cover" to ``arbitrary clopen cover" in the definition if we want. Let us now see that $\shR$ really is closely related to usual sheaf categories.

%Given a profinite space $X$, we can define a (basis for a) Grothendieck topology on $O_c(X)$ by declaring a cover to be a finite clopen cover i.e.\ given a clopen $U\subseteq X$, we say that a finite collection $\{U_i\}$ of clopen subspaces covers $U$ if $U=\bigcup U_i$. It is clear that this defines a Grothendieck topology on $O_c(X)$ and that a sheaf on this site with values in $\DModR$ is what we have defined in Definition \ref{defsh}.

First observe that we could have just considered sheaves $O(X)^{\mathrm{op}}\to\DModR$, where $O(X)$ is the poset category of open subspaces of $X$, instead of $O_c(X)^{\mathrm{op}}\to\DModR$. Indeed, there is an obvious map of Grothendieck sites $O_c(X)\to O(X)$ (given by the inclusion of clopen subsets into open subsets). This induces an equivalence between sheaves on $O(X)$ and sheaves on $O_c(X)$: as every open subset of a profinite space is a union of clopen subsets (see \cite[Theorem 1.1.12]{profinite}), the value of a sheaf $\mathcal{F}(U)$ on an open $U$ is uniquely determined by its values $\mathcal{F}(V)$ on clopen subsets $V$. One benefit, however, of using $O_c(X)$ is that the sheaf condition can be reduced to one which mentions only finite disjoint covers, rather than arbitrary covers.

Suppose we have two sheaves $(\M,X)$ and $(\Nc,Y)$. A map of Grothendieck sites from $O_c(X)$ to $O_c(Y)$ is precisely given by a continuous map $f\colon Y\to X$. Let $\shR'$ be the category whose objects are sheaves $(\M,X)$, and where a morphism $(\Theta,f)\colon(\M,X)\to(\Nc,Y)$ consists of a continuous map $f\colon Y\to X$ (i.e.\ a map of sites) and a natural transformation $\Theta\colon\M\to\Nc\circ f^{-1}$.
\[\begin{tikzcd}
O_c(X)^{\mathrm{op}} \arrow[rd, "f^{-1}"'] \arrow[rr, "\M"] & {} \arrow[d, "\Theta", Rightarrow]      & \DModR \\
                                                          & O_c(Y)^{\mathrm{op}} \arrow[ru, "\Nc"'] &       
\end{tikzcd}\]
This is the most natural way to define a category of sheaves where the space can vary. Although it seems a priori that $\shR'$ and $\shR$ are different, they turn out to be the same:

\begin{prop}\label{defsame}
The categories $\shR'$ and $\shR$ are equivalent (even isomorphic).
\begin{proof}
Define a functor $\alpha\colon\shR'\to\shR$ which is the identity on objects. Given a morphism $(\Theta,f)\colon(\M,X)\to(\Nc,Y)$ in $\shR'$, define $\alpha(\Theta,f)=(\Phi,f)\colon(\M,X)\to(\Nc,Y)$ in $\shR$, where $\Phi\colon\M\circ p_1\to\Nc\circ f^*$ is the following natural transformation. Its component $\Phi_{(U,V)}$ on $(U,V)\in O_c(X)^\mathrm{op}\times O_c(Y)^\mathrm{op}$ is the composite $$\M(U)\overset{\Theta_{U}}{\longrightarrow}\Nc(f^{-1}(U))\to\Nc(f^{-1}(U)\cap V).$$

On the other hand, define a functor $\beta\colon\shR\to\shR'$ which is also the identity on objects. Given a morphism $(\Phi,f)\colon(\M,X)\to(\Nc,Y)$ in $\shR$, define $\beta(\Phi,f)=(\Theta,f)\colon(\M,X)\to(\Nc,Y)$ in $\shR'$, where $\Theta\colon\M\to\Nc\circ f^{-1}$ is the following natural transformation. Its component on $U\in O_c(X)$ is simply $\Theta_U=\Phi_{(U,Y)}$.

It is easy to check that $\alpha$ and $\beta$ are inverse to each other.
\end{proof}
\end{prop}

From now on, we will work mostly with $\coshR$ rather than $\shR$, because that is the viewpoint of \cite{melfree} and \cite{ribesgraph}. Results about sheaves of discrete $R$-modules may be obtained by duality. By \cite[Theorem 3.3]{gareth}, the category $\mathbf{CoSh}(R)$ is equivalent to the category of ``bundles\footnote{These are called \emph{sheaves} of profinite $R$-modules in \cite{ribesgraph}, but as pointed out in \cite{gareth}, they should perhaps be called cosheaves in modern terminology, dual to sheaves of discrete $R$-modules. To avoid confusion, we will refer to them as ``bundles", which is consistent with category-theoretic language.} of profinite $R$-modules", analogous to bundles of profinite groups defined in \cite[Chapter 5]{ribesgraph}. Let us define ``bundles of $R$-modules" here and explain the equivalence briefly. Bundles should be regarded as the more algebraic viewpoint of cosheaves.

\begin{defn}[\cite{ribesgraph} Section 5.1]\label{fakesh}
Let $R$ be a profinite ring. A \emph{bundle of (profinite, right) $R$-modules} is a triple $(M,p,X)$, where $M$ and $X$ are profinite spaces and $p\colon M\sur X$ is a continuous surjection, such that each fibre $M_x=p^{-1}(x)$ has the structure of a profinite $R$-module making the following maps continuous:
\begin{eqnarray*}
M\times R\to M,&&(m,r)\mapsto mr,\\
M\times_X M\to M,&&(m,n)\mapsto m+n.
\end{eqnarray*}
A \emph{morphism of bundles} $(M,p,X)\to(N,q,Y)$ consists of continuous maps $\psi\colon M\to N$ and $f\colon X\to Y$ such that $q\psi=fp$ and such that each $\psi|_{x}\colon M_x\to N_{f(x)}$ is an $R$-module homomorphism.

Every profinite $R$-module can be canonically viewed as a bundle over the one-point space. Given a bundle $(M,p,X)$, its \emph{profinite direct sum} is a profinite $R$-module $\bigcoplus_{x\in X} M_x$ together with a morphism $M\to \bigcoplus_{X} M_x$, which is universal with respect to this property i.e.\ such that for every profinite $R$-module $N$ and every bundle morphism $M\to N$, there is a unique map (of profinite $R$-modules) $\bigcoplus_{X} M_x\to N$ making the relevant triangle commute.

(Sometimes, we will abuse notation and write $M$ for a bundle instead of $(M,p,X)$ when $X$ and $p$ are understood. In fact, we have already done this in the above definition!)
\end{defn}

\begin{prop}[\cite{ribesgraph} Proposition 5.1.2]
Let $(M,p,X)$ be a bundle of profinite $R$-modules. Then its profinite direct sum $\bigcoplus_XM_x$ exists (and can be constructed as a completion of the abstract direct sum).
\end{prop}

Remark: We are using the symbol $\coplus\,$ for the profinite direct sum as it is a completion of the usual direct sum. The symbol $\oplus$ is used in \cite{melfree}, \cite{ribesgraph} and \cite{melaspherical}, while the symbol $\boxplus$ is used in \cite{gareth} and \cite{garethrel}. Results in \cite[Chapter 5]{ribesgraph} are mostly stated for bundles of profinite groups rather than modules. The corresponding results for modules can be found in \cite[Appendix A]{garethrel}.

\begin{exmp}\label{exsh}
The following fundamental examples of bundles and profinite direct sums are given in \cite[Examples 5.1.1, 5.1.3 and 5.6.4]{ribesgraph}.
\begin{enumerate}[label=(\roman*)]
\item\label{exsh1} Suppose we have a family of profinite $R$-modules $\{M_x\}_{x\in X}$ indexed by a finite set $X$. Then we obtain a bundle $\left( \bigsqcup_XM_x,p,X\right)$, where $\bigsqcup_XM_x$ has the disjoint union topology and $p\colon\bigsqcup_XM_x\to X$ is the obvious map. The profinite direct sum $\bigcoplus_XM_x$ is the usual finite direct sum $\bigoplus_XM_x$.
\item\label{exsh2} Suppose $M$ is any profinite $R$-module and $X$ is any profinite space. Then we obtain a bundle $(M\times X,p,X)$, where $p\colon M\times X\to X$ is the natural projection from the product space. This is called the \emph{constant bundle with value $M$} and its profinite direct sum $\bigcoplus_XM\times \{x\}$ is the profinite tensor product $M\crotimes R\llbracket X\rrbracket$, where $R$ acts on $R\llbracket X\rrbracket$ on the right.
\item\label{exsh3} Suppose we have a family of profinite $R$-modules $\{M_i\}_{i\in I}$ indexed by a (possibly infinite) set $I$. Let $X=I\sqcup\{\infty\}$ be the one-point compactification of the discrete space $I$. There is a canonical profinite topology we can define on $M=\left(\bigsqcup_IM_i\right)\sqcup\{0\}$ so that $(M,p,X)$ becomes a bundle, and such that the profinite direct sum $\bigcoplus_XM_x$ is the direct \emph{product} $\prod_IM_i$ (see \cite[Examples 5.1.1(c) and 5.6.4(c)]{ribesgraph}). We will give more details of this topology in Example \ref{exproj}, since it is there that we actually need these details.

Crucially, although the profinite direct sum should be viewed as a notion of direct sums for profinite modules, it actually generalises categorical products of profinite modules, which are known to be well-behaved.
\end{enumerate}
\end{exmp}

We now describe the equivalence between $\mathbf{CoSh}(R)$ (the category of cosheaves of profinite $R$-modules) and the category of bundles of $R$-modules; more details can be found on \cite[pages 9-10]{gareth}. Given a bundle $(M,p,X)$, we obtain a cosheaf $(\M,X)$ given by $\M(U)=\bigcoplus_{u\in U}M_u$ for each clopen $U\subseteq X$. Conversely, given a cosheaf $(\M,X)$, we obtain a bundle $(M,p,X)$, where $M=\bigsqcup_{x\in X}\left(\varprojlim_{U\ni x}\M(U)\right)$ and $p\colon M\to X$ is the obvious projection map. This equivalence identifies the space of cosections $\M(U)$ ($U\subseteq X$ clopen) of a cosheaf $\M$ with the partial profinite direct sum $\bigcoplus_UM_u$ of the corresponding bundle, and identifies the fibre $M_x$ ($x\in X$) of a bundle $(M,p,X)$ with the costalk $\varprojlim_{U\ni x}\M(U)$ of the corresponding cosheaf. Henceforth, we will freely use this equivalence.

\subsection{Properties of $\coshR$}
\label{sec2.2}

We will now investigate the categorical properties of $\coshR$, illustrating the usefulness of having two perspectives on $\coshR$. Very often, we will prove a result from only one perspective. What we implicitly mean is that the result is also true in the other (equivalent) category using the equivalence just discussed.

\begin{prop}\label{propff}
The functor $\PModR\to\mathbf{CoSh}(R)$ (which sends $M$ to the cosheaf $M$ over the one-point space) is fully faithful and has a left adjoint which is given by:
\begin{itemize}
\item the ``global cosections functor" $(\M,X)\mapsto\M(X)$ from the viewpoint of cosheaves;
\item the profinite direct sum functor $\bigcoplus\,$ from the viewpoint of bundles.
\end{itemize}
Thus, $\PModR$ is a reflective subcategory of $\mathbf{CoSh}(R)$.
\begin{proof}
This follows directly from definitions.
\end{proof}
\end{prop}

\begin{prop}
The category $\mathbf{CoSh}(R)$ has all (small) limits and finite coproducts.
\begin{proof}
We use the perspective of bundles. The limit of the bundles $(M_i,p_i,X_i)$ is given by $(\lim M_i,\lim p_i,\lim X_i)$, where the latter limits are taken in $\Top$ (or equivalently in $\Pro$). This works because limits commute with limits and because the forgetful functor $\PModR\to\Pro$ preserves limits. Finite coproducts are given by $(M,p,X)\sqcup(N,q,Y)=(M\sqcup N,p\sqcup q,X\sqcup Y)$. Notice that limits are computed ``pointwise", in the sense that $(\lim M_i)_x=\lim(M_i)_{x_i}$, where $x=(x_i)\in\lim X_i\subseteq \prod X_i$. A similar statement can be said about finite coproducts.

A subtle point is that it may not be immediately obvious that $\lim p_i\colon\lim M_i\to\lim X_i$ is surjective. This is clear for products and less clear for equalisers, but it is true because bundle maps are module homomorphisms on fibres.
\end{proof}
\end{prop}

\begin{prop}[\cite{gareth} Theorem 4.6]
Profinite direct sums commute with restriction of scalars. That is, suppose $(M,p,X)$ is a bundle of $R$-modules and there is a map $S\to R$ of profinite rings (which induces a functor $\lambda\colon\coshR\to\mathbf{CoSh}(S)$). Then $\lambda(\bigcoplus_X{M_x})=\bigcoplus_X(\lambda M)_x$ as $S$-modules.
\begin{proof}
This is obvious from the perspective of cosheaves, since the profinite direct sum functor simply sends a cosheaf $(\M,X)\in\coshR$ to its global cosections space $\M(X)\in\PModR$. Now observe that the equivalence from cosheaves to bundles respects restriction of scalars.
\end{proof}
\end{prop}

\begin{prop}[\cite{ribesgraph} Proposition 5.1.6]\label{nonzero}
Let $(M,p,X)$ be a bundle. Then the canonical map $M\to \bigcoplus_{X} M_x$ is injective when restricted to each fibre $M_x$. In particular, if there is at least one non-zero fibre $M_x$, then the profinite direct sum $\bigcoplus_{X} M_x$ is also non-zero.

(Slogan: Every costalk injects into the space of global cosections.)
\end{prop}

\begin{prop}[\cite{ribesgraph} Proposition 5.1.7]\label{propil}
Profinite direct sums commute with inverse limits.
\end{prop}

We say that a bundle $(M,p,X)$ of $R$-modules is \emph{finite} if both $M$ and $X$ are finite spaces.

\begin{prop}[\cite{ribesgraph} Theorem 5.3.4]\label{propfi}
Every bundle $(M,p,X)$ is an inverse limit of finite bundles $(M_i,p_i,X_i)$, where the inverse system can be chosen such that the projection maps $(M,p,X)\to(M_i,p_i,X_i)$ are surjective.
\begin{proof}
This follows from the proof of \cite[Theorem 5.3.4]{ribesgraph}, but let us give a quick explanation here. Note that this proof is not present in the published version. For the purpose of this proof, we shall assume \cite[Theorem 5.3.4]{ribesgraph}.

Given a bundle $(M,p,X)$, we can use (the module analogue of) \cite[Theorem 5.3.4]{ribesgraph} to obtain a directed poset $I$ and an inverse system $(M_i,p_i,X_i)$ of finite bundles, where $X=\varprojlim X_i$ and $M_i=\bigsqcup_{x\in X_i} M_{i,x}$ (where $M_{i,x}$ corresponds to $A_{i,\tau}$ in \cite{ribesgraph}). We will leave it to the reader to check that $(M,p,X)=\varprojlim(M_i,p_i,X_i)$.
\end{proof}
\end{prop}

We note that the following is slightly different from the published version, though the statement in the published version is still correct.

\begin{lemma}\label{lemfac}
Let $(M,p,X)=\varprojlim_{i\in I}(M_i,p_i,X_i)$ be an inverse limit of finite bundles and let $(N,q,Y)$ be a finite bundle. Then every morphism $$(\psi,f)\colon(M,p,X)\to(N,q,Y)$$ factors through some $(M_l,p_l,X_l)$.
\begin{proof}
When $Y=*$, this was proved in \cite[Lemma 5.1.8]{ribesgraph}. For the general case, it is easy to see that there exists $j\in I$ such that $(\psi,f)$ factors through continuous maps $\psi_j\colon M_j\to N$, $f_j\colon X_j\to Y$ which are compatible with $p_j$ and $q$. The only issue is that the restrictions $\psi_j|_{x}$ might not be $R$-module homomorphisms for $x\in X_j$. But by the case $Y=*$, we can find some $k\geq j$ such that the composite $$(M,p,X)\to(N,q,Y)\xrightarrow{g}\prod_YN_y$$ factors through a bundle map $h\colon(M_k,p_k,X_k)\to\prod_YN_y$ (see the diagram at the bottom of the proof). Composing $\psi_j$ and $f_j$ with the transition maps of the inverse system gives continuous maps $\psi_k\colon M_k\to N$ and $f_k\colon X_k\to Y$.

We still cannot yet guarantee that $\psi_k$ is a homomorphism on fibres, which would follow from $g\circ (\psi_k,f_k)=h$. But both of these maps are factorisations of $g\circ(\psi,f)$ into the finite space $\prod_YN_y$, so after passing to some $l\geq k$, they become equal, which completes the proof.
\[\begin{tikzcd}
                          &                                                   & {(M,p,X)} \arrow[rd, "{(\psi,f)}"] \arrow[d] \arrow[ld] \arrow[lld] &                     &                \\
{(M_l,p_l,X_l)} \arrow[r] & {(M_k,p_k,X_k)} \arrow[rrr, bend right,"h"] \arrow[r] & {(M_j,p_j,X_j)} \arrow[r, "{(\psi_j,f_j)}",']                         & {(N,q,Y)} \arrow[r,',"g"] & \prod_YN_y
\end{tikzcd}\]
\end{proof}
\end{lemma}

\begin{corr}\label{procomp} The category $\coshR$ is the pro-completion of the category of finite bundles of $R$-modules.
\begin{proof}
In view of the last two results, it only remains to show that if a map $$(\psi,f)\colon(M,p,X)=\varprojlim(M_i,p_i,X_i)\to(N,q,Y),$$ where the $(M_i,p_i,X_i)$ and $(N,q,Y)$ are finite bundles, factors through two stages $(M_j,p_j,X_j)$ and $(M_k,p_k,X_k)$, then we can find some $l\geq j,k$ such that the two factorisations agree at stage $l$. This is obvious (since profinite spaces is the pro-completion of finite spaces) and was already used at the end of the proof of the last lemma.
\end{proof}
\end{corr}

Recall that there is a profinite tensor product $\crotimes$ on $\PModR$. It was shown in \cite[Corollary 9.1.2]{ribesgraph} that the profinite tensor product commutes with profinite direct sums. We can define a tensor product (also denoted by $\crotimes$) on $\coshR$ which extends the one on $\PModR$, but since this is not needed in the rest of the paper (except to simplify the proof of Theorem \ref{mackey}), we will develop it in the appendix.

\subsection{Key Results}
\label{sec2.3}

In this subsection, we will see that the projective dimension of a profinite direct sum is entirely determined by its fibres (Corollary \ref{impcorr}), as well as a general Mackey's decomposition formula for profinite groups (Theorem \ref{mackey}).

We mentioned at the end of the last subsection that the profinite direct sum $\coplus\,$ commutes with $\crotimes$. It turns out that profinite direct sums don't just commute with $\crotimes$, but more generally with $\Tor_i^R$ for every $i$ (\cite[Theorem 9.1.1]{ribesgraph}). The dual notion to profinite direct sums is termed \emph{continuous products} (\cite[Definition 2.4]{gareth}). Given a sheaf $(A,p,X)$ of discrete $R$-modules, we denote its continuous product (i.e.\ space of global sections) by $\oprod_XA_x$. Continuous products satisfy properties dual to those satisfied by profinite direct sums. In addition, continuous products commute with $\Ext^i_R(M,-)$, where $M$ is a profinite $R$-module (\cite[Corollary 4.7]{gareth}). In fact, we can prove slightly more. The reader should refer to \cite{profinite}, \cite{brumer} or \cite{compactmodules} for background on the homological algebra of profinite modules. We point out that the $\Ext^i_R(-,-)$ functors used below have as inputs profinite $R$-modules in the first coordinate and discrete $R$-modules in the second.

\begin{prop}\label{pext} Let $R$ be a profinite ring and $i$ be a non-negative integer.
\begin{enumerate}[label=(\roman*)]
\item\label{pext1} Suppose $(M,p,X)$ is a bundle of profinite right $R$-modules and $N$ is a profinite left $R$-module. Then $$\Tor^R_i\left(\bigcoplus_XM_x,N\right)=\bigcoplus_X\Tor^R_i(M_x,N).$$
\item\label{pext2} Suppose $(M,p,X)$ is a bundle of profinite $R$-modules and $A$ is a discrete $R$-module. Then $$\Ext^i_R\left(\bigcoplus_XM_x,A\right)=\oprod_X\Ext^i_R(M_x,A).$$
\item\label{pext3} Suppose $(A,p,X)$ is a sheaf of discrete $R$-modules and $M$ is a profinite $R$-module. Then $$\Ext^i_R\left(M,\oprod_XA_x\right)=\oprod_X\Ext^i_R(M,A_x).$$
\end{enumerate}
\begin{proof}
Parts \ref{pext1} and \ref{pext3} were proved in \cite[Theorem 9.1.1]{ribesgraph} and \cite[Corollary 4.7]{gareth} respectively. For \ref{pext2}, we use the perspective of cosheaves and sheaves. Note that $\Ext^i_R(-,A)$ induces a functor $\coshR^{\mathrm{op}}\to\mathbf{Sh}(\widehat{\Z})$ by defining $\Ext^i_R(\M,A)(U)=\Ext^i_R(\M(U),A)$. We obtain statement \ref{pext2} by taking global (co)sections.
\end{proof}
\end{prop}

\begin{corr}[Key Corollary]\label{impcorr} Let $R$ be a profinite ring.
\begin{enumerate}[label=(\roman*)]
\item\label{impcorr2} If $(M,p,X)$ is a bundle of profinite $R$-modules, then $\bigcoplus_XM_x$ has projective dimension $\leq d$ if and only if each fibre $M_x$ has projective dimension $\leq d$. In particular, any product of projective profinite modules is projective.\footnote{The author would like to thank the referee for pointing out that a profinite $R$-module is projective iff it is flat with respect to $\crotimes$ (cf.\ \cite[Proposition 3.1]{brumer}). Thus, this statement is trivially true when ``projective" is replaced with ``flat" and certainly does not require a separate proof using Proposition \ref{pext}\ref{pext1}.}
\item\label{impcorr3} If $(A,p,X)$ is a sheaf of discrete $R$-modules, then $\oprod_XA_x$ has injective dimension $\leq d$ if and only if each fibre $A_x$ has injective dimension $\leq d$. In particular, any direct sum of injective discrete modules is injective.
\end{enumerate}
(Slogan for \ref{impcorr2}: Projective dimension of space of global cosections equals maximum projective dimension of costalks.)
\begin{proof}
This was basically proved in \cite[Lemma 1.6]{melaspherical}. Part \ref{impcorr3} is the dual of \ref{impcorr2}. Moreover, for the ``in particular" claims, note that $\prod_IM_i=\bigcoplus_{I^*}M_i$, where $I^*$ is the one-point compactification of the discrete index set $I$ (see Example \ref{exsh}\ref{exsh3}). The ``in particular" claims are well-known (see \cite[Corollary 3.3]{brumer}).

For \ref{impcorr2}, note that a profinite $R$-module $N$ has projective dimension $\leq d$ if and only if $\Ext^{d+1}_R(N,F)=0$ for all finite $R$-modules $F$. The case $d=0$ follows from \cite[Lemma 5.4.1]{profinite}, and the general case follows by dimension shifting. Now use Proposition \ref{pext}\ref{pext2}. Recall that by (the dual of) Proposition \ref{nonzero}, a continuous product is zero if and only if every stalk is zero.
\end{proof}
\end{corr}

\begin{rmk}
\begin{enumerate}[label=(\roman*)]
\item Let $(M,p,X)$ be a bundle. In the language of cosheaves, it is clear from Corollary \ref{impcorr}\ref{impcorr2} that every costalk $M_x$ ($x\in X$) has projective dimension $\leq d$ if and only if \emph{every} cosections space $\bigcoplus_UM_u$ ($U\subseteq X$ clopen) has projective dimension $\leq d$. Thus, the property of ``having projective dimension $\leq d$" transfers between costalks and cosections spaces (i.e.\ it is in this sense a (co)local property!). This is analogous to the fact that for abstract modules, a direct sum has projective dimension $\leq d$ if and only if every summand does.
\item It is natural to ask if the property of ``having projective dimension $\leq d$" also transfers between costalks and the entire cosheaf $M$. For this question to make sense, fix a profinite space $X$ and consider the category $\mathbf{CoSh}(R,X)$ of usual cosheaves on $X$ with coefficients in $\PModR$. (Warning: This is \textbf{not} the full subcategory of $\coshR$ where the profinite space happens to be $X$, because we do not allow any non-identity endomorphisms of $X$ in $\mathbf{CoSh}(R,X)$.)

It is easy to check that $\mathbf{CoSh}(R,X)$ is an abelian subcategory of the precosheaf category $[O_c(X),\PModR]$ closed under all limits and colimits. We can ask whether it is true that for any cosheaf $\M\in\mathbf{CoSh}(R,X)$, every costalk $\M_x\in\PModR$ is projective if and only if $\M\in\mathbf{CoSh}(R,X)$ is itself projective. The answer is no, as the next example shows. Thus, Corollary \ref{impcorr}\ref{impcorr2} ``skips" the level of cosheaves and directly relates the projective dimensions of costalks in $\PModR$ with the projective dimensions of cosections spaces in $\PModR$.
\end{enumerate}
\end{rmk}

\begin{exmp}\label{exproj}
Let $X=\Z^+\sqcup\{\infty\}$ be the one-point compactification of $\Z^+$ (the discrete space of positive integers) and suppose $P\in\PModR$ is any non-zero projective module (say $P=R$). Consider two cosheaves $\M$ and $\Nc$ on $X$ as follows. The cosheaf $\M$ is the constant cosheaf with value $P$. Explicitly, from the perspective of bundles, it is given by $(P\times X,p,X)$, where $p\colon P\times X\to X$ is the natural projection from the product space (see Example \ref{exsh}\ref{exsh2}). Every fibre (i.e.\ costalk) of $\M$ is $P$, which is projective.

On the other hand, $\Nc$ is the \emph{cosheaf with value $P$ converging to $0$}. Explicitly, from the perspective of bundles, it is given by $((P\times \Z^+)\sqcup\{0\},q,X)$, with the following topology on $(P\times \Z^+)\sqcup\{0\}$. A subset $U\subseteq (P\times \Z^+)\sqcup\{0\}$ not containing $0$ is open if and only if it is open in the product $P\times \Z^+$, whilst subsets of the form $(P\times Y)\sqcup\{0\}$, where $Y\subseteq\Z^+$ is cofinite, form a base of open sets around $0$. The map $q$ projects $P\times \Z^+$ onto $\Z^+$ and sends $0$ to $\infty$ (cf. Example \ref{exsh}\ref{exsh3}). Every fibre of $\Nc$ is either $P$ or $0$, so is projective.

There is a natural epimorphism $\M\sur\Nc$ in $\mathbf{CoSh}(R,X)$ which does not split. Indeed, the only possible candidate for such a splitting $\psi\colon (P\times \Z^+)\sqcup\{0\}\to P\times X$ has to be the identity on $P\times \Z^+$ and send $0$ to $(0,\infty)\in P\times X$. However, $\psi$ is not continuous: if $Q\subsetneq P$ is any open subset containing 0, then $\psi^{-1}(Q\times X)=(Q\times \Z^+)\sqcup\{0\}$ is not open in $(P\times \Z^+)\sqcup\{0\}$. Hence, $\Nc\in\mathbf{CoSh}(R,X)$ fails to be projective, even though all of its fibres are projective.
\end{exmp}

Next, we shall prove the profinite Mackey's Formula. Let $R$ be a ring, $G$ be a group, and $X$ be a (right) $G$-set. Then we can always decompose the $R[G]$-module $R[X]$ as $R[X]=\bigoplus_{x\in X/G}R[xG]$. If $R$, $G$, $X$ are all profinite, then a direct analogue of this works if $X/G$ is finite. In general, we can still get such a decomposition, but we need to use the profinite direct sum instead of the usual direct sum.

\begin{prop}[\cite{melaspherical} Corollary 1.11/\cite{garethrel} Proposition A.14]\label{meldec}
Let $R$ be a profinite ring, $G$ a profinite group, and $X$ a profinite $G$-space. Then there is an isomorphism of profinite $\RG$-modules $$R\llbracket X\rrbracket=\bigcoplus_{x\in X/G}R\llbracket xG\rrbracket.$$
\end{prop}

The same method shows that we get a Mackey's decomposition formula for profinite groups if we use the profinite direct sum. This was essentially done in \cite[Section 3]{double} but without the language of cosheaves. Although most of our modules so far have been one-sided, there is a completely analogous way to define bundles of $R$-$S$-bimodules and their profinite direct sums. See also \cite[Chapter 6]{profinite} for some background on induction and restriction.

\begin{thm}[Profinite Mackey's Formula]\label{mackey}
Let $G$ be a profinite group and $H,K$ closed subgroups (which are not necessarily open). If $M$ is a profinite right $R\llbracket H\rrbracket$-module, then there is an isomorphism of profinite right $R\llbracket K\rrbracket$-modules $$\Res^G_K\Ind^G_HM=\bigcoplus_{g\in H\backslash G/K}\Ind^K_{K\cap g^{-1}H g}\Res^{g^{-1}H g}_{K\cap g^{-1}H g}Mg.$$
\begin{proof}
There is an $R\llbracket H\rrbracket$-$R\llbracket K\rrbracket$-bimodule isomorphism $$\RG=\bigcoplus_{g\in H\backslash G/K}R\llbracket HgK\rrbracket,$$ which is easy to prove by following the argument of \cite[Proposition A.14]{garethrel}. We then apply the functor $M\,\,\widehat{\otimes}_{R\llbracket H\rrbracket}\,\,(-)$ (from the category of $R\llbracket H\rrbracket$-$R\llbracket K\rrbracket$-bimodules to the category of right $R\llbracket K\rrbracket$-modules) to the above. Note that $M\,\,\widehat{\otimes}_{R\llbracket H\rrbracket}\,\,(-)$ commutes with profinite direct sums of bimodules: one may prove this directly just like in Proposition \ref{pext}\ref{pext1}, but in Remark \ref{rm}\ref{rm4} of the appendix, we will provide a formal reason why this follows from the corresponding result for one-sided modules (Proposition \ref{pext}\ref{pext1}).

It remains to show that $M\,\,\widehat{\otimes}_{R\llbracket H\rrbracket}\,\,R\llbracket HgK\rrbracket=\Ind^K_{K\cap g^{-1}H g}\Res^{g^{-1}H g}_{K\cap g^{-1}H g}Mg$ as right $R\llbracket K\rrbracket$-modules, which would follow from the isomorphism $R\llbracket HgK\rrbracket=R\llbracket Hg\rrbracket\,\,\widehat{\otimes}_{R\llbracket K\cap g^{-1}Hg\rrbracket}\,\,R\llbracket K\rrbracket$. To prove the last isomorphism, write $G=\varprojlim_iG_i$ ($G_i$ finite), with surjective projection maps $\varphi_i\colon G\sur G_i$. Let $H_i=\varphi_i(H)$, $K_i=\varphi_i(K)$ and $g_i=\varphi_i(g)$. Then $H=\varprojlim H_i$, $K=\varprojlim K_i$ and $HgK=\varprojlim\varphi_i(HgK)=\varprojlim H_ig_iK_i$. Moreover, we have $$K\cap g^{-1}Hg=\varprojlim\varphi_i(K)\cap\varprojlim\varphi_i(g^{-1}Hg)=\varprojlim(\varphi_i(K)\cap\varphi_i(g^{-1}Hg))=\varprojlim (K_i\cap g_i^{-1}H_ig_i).$$ Combining all of these, we see that
\begin{eqnarray*}
R\llbracket Hg\rrbracket\,\,\widehat{\otimes}_{R\llbracket K\cap g^{-1}Hg\rrbracket}\,\,R\llbracket K\rrbracket&=&\varprojlim \left(R[H_ig_i]\otimes_{R[ K_i\cap g_i^{-1}H_ig_i]}R[ K_i]\right)\\
&=&\varprojlim R[H_ig_iK_i]\\
&=&R\llbracket HgK\rrbracket,
\end{eqnarray*}
as required.
\end{proof}
\end{thm}

Remark: If $H$ or $K$ is open, so that the quotient $H\backslash G/K$ is finite, then $\coplus\,$ agrees with $\oplus$ (Example \ref{exsh}\ref{exsh1}) and we get the standard Mackey's decomposition formula (see \cite[Corollary 2.2]{double}).

\section{Profinite Groups of Type $\Phi_R$}
\label{sec3}

In this section, we will apply results from Section \ref{sec2.3} to profinite groups of type $\Phi_R$. First, let us explain why one might care about these groups.

\subsection{Motivations and Definition}
\label{sec3.1}

Groups of type $\Phi$ were first introduced in \cite{talephi} as a potential algebraic characterisation of those groups $G$ which admit a finite-dimensional model for $\underline{E}G$, the classifying space for finite subgroups of $G$. We say that an abstract group $G$ is \emph{of type $\Phi_R$ (over a ring $R$)} if an $R[G]$-module has finite projective dimension whenever its restriction to every finite subgroup of $G$ has finite projective dimension. It was later shown in \cite[Theorem 3.10]{infinitestable} that over many rings, these groups have well-behaved ``stable module categories", which make their representation theory easier to study. There is an obvious way to adapt the definition of type $\Phi_R$ to the profinite setting.

\textbf{Throughout Section \ref{sec3}, $R$ will be a commutative profinite ring of finite global dimension $n$} (i.e.\ every profinite $R$-module has projective dimension at most $n$). Important examples include $\F_p$ and $\Z_p$.

\begin{defn}\label{phidef}
Let $G$ be a profinite group. We say that $G$ is \emph{of type $\Phi_R$} if a profinite $\RG$-module has finite projective dimension whenever its restriction to every finite subgroup of $G$ has finite projective dimension.
\end{defn}

In the definition above, we could have changed the last part to ``whenever its restriction to every finite subgroup of $G$ has projective dimension $\leq n$" (recall that $n$ is the global dimension of $R$). This is because of the following lemma.

\begin{lemma}\label{peterlem}
Suppose $G$ is a finite group and $M$ is a profinite $\RG$-module of finite projective dimension. Then $M$ has projective dimension $\leq n$.
\begin{proof}
This is the same argument as \cite[Lemma 2.3]{infinitestable}, but we include it for completeness. Suppose instead that $M$ has finite $\RG$-projective dimension $m>n$. Consider the $\RG$-short exact sequence $P\inj Q\sur K$, where $P$ and $Q$ are projective and $K$ is a (non-projective) $(m-1)$th syzygy of $M$. Then $K$ is projective over $R$, so this sequence is $R$-split, but $P$ is projective over $\RG$ and $G$ is finite, so $P$ is injective relative to the trivial subgroup. Hence, this sequence is actually $\RG$-split and $K$ is $\RG$-projective, a contradiction.
\end{proof}
\end{lemma}

\subsection{Properties of type $\Phi_R$}
\label{sec3.2}

We now study the properties of the class of profinite groups of type $\Phi_R$. All of these properties have abstract counterparts, whose proofs often use infinite direct sums in seemingly essential ways. Fortunately, the profinite direct sum is enough to play the role of the abstract direct sum in the following results. In the next subsection, we will show that the class of profinite groups of type $\Phi_R$ is fairly large.

\begin{lemma}\label{fdlem}
Suppose $G$ is a profinite group of type $\Phi_R$. Then there exists some constant $c$ such that every profinite $\RG$-module of finite projective dimension has projective dimension at most $c$.
\begin{proof}
The idea is similar to \cite[Lemma 2.4]{infinitestable}, but we need to use products instead of direct sums. Recall that products are just special cases of profinite direct sums (Example \ref{exsh}\ref{exsh3}). Suppose instead that for each positive integer $i$, there is a module $M_i$ of finite projective dimension at least $i$. Then the product $M=\prod_iM_i$ has projective dimension $\leq n$ when restricted to any finite subgroup of $G$, but $M$ does not have finite projective dimension over $\RG$, a contradiction.
\end{proof}
\end{lemma}

\begin{defn}
The smallest such $c$ from the above lemma is called the \emph{finitistic dimension} of $\RG$ and denoted by $\Fd\RG$.
\end{defn}

Remark: In general, for a profinite ring $S$, we define its \emph{finitistic dimension} $\Fd S$ to be the supremum of the projective dimensions of all $S$-modules with finite projective dimension. Lemma \ref{fdlem} shows that if a profinite group $G$ is of type $\Phi_R$, then $\Fd\RG<\infty$.

The following is the profinite analogue of \cite[Proposition 2.3(i)]{talephi} and is a first application (to the author's knowledge) of the general Mackey's decomposition formula for profinite groups (Theorem \ref{mackey}).

\begin{prop}\label{phisub}
If a profinite group $G$ is of type $\Phi_R$ then so is every closed subgroup $H$ of $G$. Moreover, $\Fd R\llbracket H\rrbracket\leq\Fd\RG$.
\begin{proof}
Let $\Fd\RG=c$. Suppose $M$ is an $R\llbracket H\rrbracket$-module whose restriction to every finite subgroup of $H$ has projective dimension $\leq n$. We want to show that $M$ has $R\llbracket H\rrbracket$-projective dimension $\leq c$. If $K$ is any closed subgroup of $G$, then by Mackey's Formula (Theorem \ref{mackey}), we have $$\Res^G_K\Ind^G_HM=\bigcoplus_{g\in H\backslash G/K}\Ind^K_{K\cap g^{-1}H g}\Res^{g^{-1}H g}_{K\cap g^{-1}H g}Mg.$$ Importantly, if $K$ is finite, then $\Res^{g^{-1}H g}_{K\cap g^{-1}H g}Mg$ has projective dimension $\leq n$ by assumption, so the entire right hand side has projective dimension $\leq n$ by Corollary \ref{impcorr}\ref{impcorr2}. As $G$ is of type $\Phi_R$, we conclude that $\Ind^G_HM$ has projective dimension $\leq c$ over $\RG$, so $\Res^G_H\Ind^G_HM$ has projective dimension $\leq c$ over $R\llbracket H\rrbracket$. We used the fact that restriction and induction preserve projective resolutions (see \cite[Corollary 5.7.2 and Theorem 6.10.8]{profinite}).

We use Mackey's Formula again with $K=H$, which reads  $$\Res^G_H\Ind^G_HM=\bigcoplus_{g\in H\backslash G/H}\Ind^H_{H\cap g^{-1}H g}\Res^{g^{-1}H g}_{H\cap g^{-1}H g}Mg.$$ Importantly, taking $g=1$ shows that the $R\llbracket H\rrbracket$-module $M$ appears as a fibre on the right hand side. By Corollary \ref{impcorr}\ref{impcorr2}, $M$ has $R\llbracket H\rrbracket$-projective dimension $\leq c$, as required.
\end{proof}
\end{prop}

Remark: In the abstract case, Mackey's decomposition formula shows in particular that a $H$-module $M$ is a direct summand of $\Res^G_H\Ind^G_HM$. In the profinite case, we may not be able to make this conclusion, but what we can say is that $M$ appears as a fibre in some profinite direct sum decomposition of $\Res^G_H\Ind^G_HM$. This is sufficient for the proof above.

\subsection{Examples of type $\Phi_R$}
\label{sec3.3}

We shall see in this subsection that for some profinite rings $R$ (such as $\F_p$), the class of profinite groups of finite $vcd_R$ is strictly contained in the class of profinite groups of type $\Phi_R$. The example we provide for the ``strictly" part will use some profinite Bass-Serre Theory, for which an excellent reference is \cite{ribesgraph}.

The published version of Corollary \ref{corrtype}\ref{corra} is slightly wrong. To fix it, we will need the following proposition, which is slightly stronger than the published version.

\begin{prop}\label{typepro}
Let $G$ be a profinite group. Suppose the trivial module $R$ has a finite length $\RG$-module resolution $0\to N_l\to\cdots\to N_0\sur R$ (of length $l$) such that:
\begin{enumerate}[label=(\roman*)]
\item Each $N_i$ is a summand of a module $R\llbracket U_i\rrbracket$, where $U_i$ is a profinite $G$-space, and the $\RG$-action on $R\llbracket U_i\rrbracket$ is induced by the $G$-action on $U_i$, and
\item The $G$-stabiliser $G_u$ of each $u\in U_i$ is finite.
\end{enumerate}
Then $G$ is of type $\Phi_R$ with $\Fd\RG\leq n+l$.
\begin{proof}
By Proposition \ref{meldec}, we can write each $R\llbracket U_i\rrbracket$ as $R\llbracket U_i\rrbracket=\bigcoplus_{u\in U_i/G}R\llbracket G_u\backslash G\rrbracket$, where each $G_u$ is finite by assumption. Suppose $M$ is an $\RG$-module whose restriction to every finite subgroup of $G$ has projective dimension $\leq n$. Note that each term in the given sequence $N_*\sur R$ is $R$-projective and hence the sequence is $R$-split. Applying the functor $M\crotimes(-)$ to the $R$-split resolution $N_*\sur R$ gives a finite length $\RG$-resolution $M\crotimes N_*\sur M$ (with $G$ acting diagonally), where each term (except $M$) is a summand of \begin{eqnarray*}M\crotimes\left(\bigcoplus_{u\in U_i/G}R\llbracket G_u\backslash G\rrbracket\right)&=&\bigcoplus_{{u\in U_i/G}}(M\crotimes R\llbracket G_u\backslash G\rrbracket)\\&=&\bigcoplus_{{u\in U_i/G}}(M\,\,\widehat{\otimes}_{R\llbracket G_u\rrbracket}\,\,R\llbracket G\rrbracket)\\&=&\bigcoplus_{{u\in U_i/G}}\Ind^G_{G_u}\Res^G_{G_u}M,\end{eqnarray*} which has projective dimension $\leq n$ over $\RG$ by Corollary \ref{impcorr}\ref{impcorr2}. The first isomorphism above follows from Proposition \ref{pext}\ref{pext1} and the second follows from \cite[Proposition 5.8.1]{profinite}. By dimension shifting, we have $\Ext_{\RG}^{n+l+1}(M,-)=\Ext_{\RG}^{n+1}(M\crotimes N_l,-)=0$, so $M$ has projective dimension $\leq n+l$ over $\RG$.
\end{proof}
\end{prop}

Recall that given a profinite group $G$, its \emph{cohomological dimension over $R$}, or \emph{$cd_R$}, is the $\RG$-projective dimension of the trivial module $R$, whilst its \emph{virtual cohomological dimension over $R$}, or \emph{$vcd_R$}, is the cohomological dimension of any open subgroup with finite $cd_R$ (or $\infty$ if no such subgroup exists).

The following is the profinite analogue of \cite[Corollary 2.6]{infinitestable}.

\begin{corr}\label{corrtype}
Let $G$ be a profinite group, $\pi$ a non-empty set of prime numbers and $p$ a fixed prime. Suppose one of the following holds:
\begin{enumerate}[label=(\alph*)]
\item\label{corra} $G$ has finite $vcd_R$, where $R$ is a commutative profinite ring (of finite global dimension) satisfying the assumptions of \cite[Section 2.1]{permcom} (this includes $R=\F_p$ and $R=\Z_p$).
\item\label{corrb} $G$ acts on a $\pi$-tree $T$ with finite (vertex and edge) stabilisers, and $R$ is a commutative pro-$\pi$ ring (of finite global dimension).
\end{enumerate}
Then $G$ is of type $\Phi_R$.
\begin{proof}
\begin{enumerate}[label=(\alph*)]
\item \cite[Corollary 6.17]{permcom} gives an $\RG$-resolution of $R$ satisfying the hypotheses of Proposition \ref{typepro}.
%Let $H$ be an open subgroup of $G$ of finite $cd_p$. Take a finite $\FH$-projective resolution $Q_*\sur R$. Since $\FH$ is a local ring, projective modules are free (see \cite[Lemmas 1.4 and 1.9]{melaspherical}). By \cite[Proposition 10.2.4, Corollary 10.2.6, Theorem 10.2.8]{ribesgraph}, we can take the complete tensor product of $Q_*$ with itself to get a finite length $\RG$-resolution $N_*\sur R$ which satisfies the hypotheses of Proposition \ref{typepro}.
\item Let $\Z_{\pi}=\prod_{q\in\pi}\Z_q$. By definition (\cite[Section 2.4]{ribesgraph}), there is a split exact sequence of $\Z_{\pi}$-modules  $$0\to \Z_{\pi}\llbracket (E^*,*)\rrbracket\to \Z_{\pi}\llbracket V\rrbracket\to \Z_{\pi}\to0,$$ where $V$ is the vertex set of $T$, $E$ is the edge set, and $E^*=T/V$. Applying $R\,\,\widehat{\otimes}_{\,\Z_{\pi}}\,(-)$ to the above sequence gives an exact sequence of $R$-modules $$0\to R\llbracket (E^*,*)\rrbracket\to R\llbracket V\rrbracket\to R\to0.$$ The action of $G$ on $T$ naturally makes the above modules into $\RG$-modules, and the stabilisers are finite by assumption. We note that the module $R\llbracket (E^*,*)\rrbracket$ is defined over a \emph{pointed} profinite space, but the argument of Proposition \ref{typepro} still works for pointed spaces. (In fact, \cite[Corollary 1.11]{melaspherical} is the pointed version of Proposition \ref{meldec}.)
\end{enumerate}
\end{proof}
\end{corr}

Remark: In part \ref{corrb} of the proof above, we used the existence of a two-term resolution of $R$. There is of course nothing special about ``two" here, so $G$ is of type $\Phi_R$ if it ``acts on a finite-dimensional contractible $\pi$-CW complex with finite stabilisers", in the sense that there is a longer exact sequence of $R$-modules similar to what we had above.

%In the abstract case, we have easy examples of groups which act on trees with finite stabilisers (and hence of type $\Phi_p$), but having infinite $vcd_p$. For example, we can take $G$ to be the countably infinite direct sum of $\Z/p$ (see \cite[page 5]{summernotes} for this example, or \cite[Example 3]{ikenaga} for the case of $\Z$ instead of $\F_p$). However, in the pro-$p$ case, an analogue of this has no chance of working. Indeed, by \cite[Theorem 4.2.11]{ribesgraph}, if an infinite direct product of $\Z/p$ acts on a $p$-tree, then there is necessarily a globally fixed vertex (which cannot have finite stabiliser). The same conclusion holds for groups such as $\Z/p\times\Z/p^2\times\Z/p^3\times\cdots$, which have torsion elements of unbounded order.

Let us now show that for some $R$, there are profinite groups of type $\Phi_R$ which have infinite $vcd_R$. Note that if $G$ is a torsion-free profinite group, then it is of type $\Phi_R$ if and only if it has finite $cd_R$, so we should only search for profinite groups with non-trivial torsion.

First, let us recall some basic profinite Bass-Serre Theory from \cite{ribesgraph}. Let $\C$ be an extension-closed pseudovariety of finite groups, that is, $\C$ is a (non-empty) class of finite groups which is closed under subgroups, images, finite products and extensions. Let $\pi=\pi(\C)$ be the set of primes which divide the order of some group in $\C$. For example, $\C$ might be the class of all finite $p$-groups, in which case $\pi(\C)=\{p\}$.

Given a connected profinite graph $\Gamma$ and a graph of pro-$\C$ groups $(\mathcal{G},\varpi,\Gamma)$ over $\Gamma$, we obtain its fundamental pro-$\C$ group $\Pi=\Pi^\C_1(\mathcal{G},\Gamma)$, the $\C$-standard graph $S=S^\C(\mathcal{G},\Gamma)$, and a canonical action of $\Pi$ on $S$ (see \cite[Sections 6.1, 6.2 and 6.3]{ribesgraph}). Moreover, the stabilisers of this action are finite if all the fibres $\mathcal{G}(m)$ of $\mathcal{G}$ are finite (\cite[Lemma 6.3.2]{ribesgraph}), and $S$ is a $\C$-tree or equivalently a $\pi$-tree (\cite[Corollary 6.3.6]{ribesgraph}). Thus, if the fibres of $\mathcal{G}$ are finite, then by Corollary \ref{corrtype}\ref{corrb}, $\Pi$ is of type $\Phi_R$, where $R$ is any commutative pro-$\pi$ ring (of finite global dimension).

An example of this is given in \cite[Example 6.2.3(b)]{ribesgraph}. Indeed, if $G=\coprod_X G_x$ is any free pro-$\C$ product of finite $\C$-groups $G_x$ (continuously indexed by the profinite space $X$), then $G=\Pi^\C_1(\mathcal{G},T)$ is the fundamental group over some $\C$-tree $T$. Furthermore, each fibre of $\mathcal{G}$ is either some $G_x$ or trivial, and so finite. Hence, we have shown the following.

\begin{prop}
Any free pro-$\C$ product of finite groups is of type $\Phi_R$, where $R$ is a commutative pro-$\pi(\C)$ ring (of finite global dimension).
\end{prop}

Remark: In Section \ref{sec2}, we defined a bundle of profinite $R$-modules and its profinite direct sum. There is a completely analogous definition for a bundle of pro-$\C$ groups and its free pro-$\C$ product (\cite[Section 5.1]{ribesgraph}), which is what we used above.

Fix a prime $p$. When $R=\F_p$, we write $vcd_p$ for $vcd_{\F_p}$ and write $\Phi_p$ for $\Phi_{\F_p}$. Next, we construct a free pro-$p$ product of finite groups which does not have finite $vcd_p$. Let $$G=\Z/p\coprod\Z/p^2\coprod\Z/p^3\coprod\cdots$$ be the free pro-$p$ product over the one-point compactification $X=\Z^+\sqcup\{\infty\}$ (see \cite[Example 5.1.1(c)]{ribesgraph}). The fact that $G$ does not have finite $vcd_p$ follows from the following simple observation:

\begin{lemma}
If a profinite group $G$ has finite $p$-subgroups of arbitrarily large orders, then it has infinite $vcd_p$.
\begin{proof}
Suppose instead that $G$ has finite $vcd_p$, with an open subgroup $H$ of finite $cd_p$. Then the normal core $H_G$ of $H$ is also open and has finite $cd_p$. Take a finite $p$-subgroup $K\leq G$ of order larger than $|G:H_G|$. Then $K\cap H_G\neq1$, so $H_G$ is not $p$-torsion-free, a contradiction to having finite $cd_p$.
\end{proof}
\end{lemma}

Hence, we finally have the following.

\begin{thm}\label{thmeg}
The class of profinite groups of finite $vcd_p$ is strictly contained in the class of profinite groups of type $\Phi_p$. Any free pro-$p$ product of finite groups with arbitrarily large torsion is of type $\Phi_p$ but not of finite $vcd_p$.
\end{thm}

\appendix
\section{Tensor Products of Cosheaves}
\label{appendix}

Let $R$ be a profinite ring. Recall that there is a profinite tensor product $\crotimes$ on $\PModR$ (see \cite[Section 5.5]{profinite}). In this appendix, we will define a tensor product $\crotimes$ on $\coshR$ which extends the one on $\PModR$. We have kept the following results separate from the main paper so that they do not interrupt its flow.

Let us first define the tensor product from the viewpoint of bundles, which more closely resembles the usual tensor product $\crotimes$ on $\PModR$. Given a finite bundle of right $R$-modules $(M,p,X)$ and a finite bundle of left $R$-modules $(N,q,Y)$, we define their \emph{tensor product}, denoted by $(M,p,X)\otimes_R(N,q,Y)\in\mathbf{CoSh}(\widehat{\Z})$, or simply $M\otimes_R N$, as follows. The second profinite space of the tensor product is $X\times Y$ and the first profinite space is $\bigsqcup_{(x,y)\in X\times Y}{M_x\otimes_R N_y}$, with the obvious projection map. This is a finite bundle of abelian groups with fibres given by $(M\otimes_R N)_{(x,y)}=M_x\otimes_R N_y$.

Given general bundles of right $R$-modules $(M,p,X)=\varprojlim_{i}(M_i,p_i,X_i)$ and left $R$-modules $(N,q,Y)=\varprojlim_{j}(N_j,q_j,Y_j)$ which are written as inverse limits of finite bundles, we define their \emph{tensor product} as $$(M,p,X)\crotimes(N,q,Y)=\varprojlim_{i,j}((M_i,p_i,X_i)\otimes_R(N_j,q_j,Y_j)),$$ which is a bundle of profinite abelian groups i.e.\ an object of $\mathbf{CoSh}(\widehat{\Z})$. Note that the second profinite space of the tensor product $(M,p,X)\crotimes (N,q,Y)$ is $X\times Y$. To describe the fibres, given $(x,y)\in X\times Y$, let us write $x=(x_i)\in\varprojlim X_i$ and $y=(y_j)\in\varprojlim Y_j$. The fibres of the tensor product are then given by $(M\crotimes N)_{(x,y)}=\varprojlim ((M_i)_{x_i}\otimes_R(N_j)_{y_j})=M_x\crotimes N_y$.

To see that the tensor product of two bundles is (up to isomorphism) independent of the way they are written as inverse limits, we need the following universal property. Given bundles of right $R$-modules $(M,p,X)$ and left $R$-modules $(N,q,Y)$, as well as a bundle of profinite abelian groups $(K,r,Z)$, an \emph{$R$-middle linear map} from $(M,p,X)\times(N,q,Y)$ to $(K,r,Z)$ consists of continuous maps $\psi\colon M\times N\to K$ and $f\colon X\times Y\to Z$ such that $r\psi=f(p\times q)$ and such that for each $(x,y)\in X\times Y$, the restriction $\psi|_{(x,y)}\colon M_x\times N_y\to K_{f(x,y)}$ is an $R$-middle linear map i.e.\ $\psi|_{(x,y)}$ is bilinear and satisfies $\psi|_{(x,y)}(mr,n)=\psi|_{(x,y)}(m,rn)$. We can now state the universal property of the tensor product, which is exactly what one might expect. We note that the proof in the published version is incorrect, so let us correct it below.

\begin{prop}\label{corrected}
Let $(M,p,X)$ and $(N,q,Y)$ be bundles of appropriately sided $R$-modules. Then their tensor product $(M,p,X)\crotimes(N,q,Y)\in\coshZ$ satisfies the following universal property:

There is an $R$-middle linear map $$(M,p,X)\times(N,q,Y)\to(M,p,X)\crotimes(N,q,Y),$$ such that for any $(K,r,Z)\in\coshZ$ and any $R$-middle linear map $(M,p,X)\times(N,q,Y)\to(K,r,Z)$, there exists a unique bundle map $(M,p,X)\crotimes(N,q,Y)\to(K,r,Z)$ in $\coshZ$ making the relevant triangle commute.
\begin{proof}
The middle linear map $(M,p,X)\times(N,q,Y)\to(M,p,X)\crotimes(N,q,Y)$ is the inverse limit of the middle linear maps $$(M_i,p_i,X_i)\times(N_j,q_j,Y_j)\to(M_i,p_i,X_i)\otimes_R(N_j,q_j,Y_j).$$ Since every bundle can be written as an inverse limit of finite bundles, it suffices to verify the universal property for a finite bundle of abelian groups $(K,r,Z)$. In this case, we can use the exact same technique as in the proof of Lemma \ref{lemfac}, including the middle linear analogue of \cite[Lemma 5.1.8]{ribesgraph}, to show that a middle linear map $(M,p,X)\times(N,q,Y)\to(K,r,Z)$ factors through some middle linear map $(M_i,p_i,X_i)\times(N_j,q_j,Y_j)\to(K,r,Z)$. As all spaces are now finite, it is clear that we have an induced map $(M_i,p_i,X_i)\otimes_R(N_j,q_j,Y_j)\to(K,r,Z)$, which in turn gives a bundle map $(M,p,X)\crotimes(N,q,Y)\to(K,r,Z)$ having the required property. Uniqueness of this map is clear.
\end{proof}
\end{prop}

We shall now describe the tensor product $\crotimes$ from the viewpoint of cosheaves. Let $X$ and $Y$ be profinite spaces. First note that every clopen subspace $U$ of $X\times Y$ can be written as a finite disjoint union $U=\bigsqcup_{i=1}^nX_i\times Y_i$, where $X_i\subseteq X$ and $Y_i\subseteq Y$ are clopen. Given two appropriately sided cosheaves $(\M,X)$ and $(\Nc,Y)$, we define $(\M\crotimes \Nc)(U)=\bigoplus_{i=1}^n(\M(X_i)\crotimes\Nc(Y_i))$. This indeed defines a cosheaf of profinite abelian groups (so is in particular functorial) and is independent of the decomposition of $U$ since we can pass to a common refinement of two given decompositions. It is not difficult to see that this tensor product $\crotimes$ on cosheaves agrees with the one we just defined on bundles by using the following proposition:

\begin{prop}[\cite{ribesgraph} Corollary 9.1.2]\label{tensorcomm}
The profinite tensor product $\crotimes$ on $\PModR$ commutes with profinite direct sums, that is, for $(M,p,X)$ and $(N,q,Y)$ appropriately sided bundles of $R$-modules, we have an isomorphism of profinite abelian groups $$\left(\bigcoplus_X M_x\right)\crotimes\left(\bigcoplus_YN_y\right)=\bigcoplus_{X\times Y}(M_x\crotimes N_y).$$
\end{prop}

\begin{corr}
The two tensor products $\crotimes$ on $\coshR$ we have just defined agree under the equivalence of cosheaves and bundles.
\begin{proof}
Recall that the equivalence $F$ from bundles to cosheaves sends a bundle $(M,p,X)$ to the cosheaf $(F(M),X)$ defined by $F(M)(V)=\bigcoplus_VM_v$, where $V\subseteq X$ is clopen. Let $(M,p,X)$ and $(N,q,Y)$ be bundles and write a clopen subspace $U\subseteq X\times Y$ as a finite disjoint union $U=\bigsqcup X_i\times Y_i$, where the $X_i$ and $Y_i$ are clopen. Then
\begin{eqnarray*}
F(M\crotimes N)(U)&=&\bigcoplus_{u\in U}(M\crotimes N)_u\\
&=&\bigoplus_i\bigcoplus_{(x,y)\in X_i\times Y_i}(M\crotimes N)_{(x,y)}\\
&=&\bigoplus_i\bigcoplus_{(x,y)\in X_i\times Y_i}(M_x\crotimes N_y)\\
&=&\bigoplus_i\left(\left(\bigcoplus_{x\in X_i}M_x\right)\crotimes\left(\bigcoplus_{y\in Y_i}N_y\right)\right)\\
&=&\bigoplus_i(F(M)(X_i)\crotimes F(N)(Y_i))\\
&=&(F(M)\crotimes F(N))(U).
\end{eqnarray*}
\end{proof}
\end{corr}

\begin{rmk}\label{rm}
\begin{enumerate}[label=(\roman*)]
\item\label{rm1} Suppose $R$ is commutative, so that $M\crotimes N$ is a cosheaf of profinite $R$-modules (rather than just profinite abelian groups). Then both $\PModR$ and $\coshR$ become symmetric monoidal categories under the tensor product $\crotimes$ with tensor unit $R$. Moreover, the adjunction of Proposition \ref{propff} is a strong monoidal adjunction, by Proposition \ref{tensorcomm}. The reader should refer to \cite{maccat} for the necessary background on monoidal categories.
\item\label{rm2} More generally, suppose $R$ is still a commutative profinite ring and let $G$ be a profinite group. Consider profinite modules over the completed group algebra $\RG$. This is a symmetric monoidal category under the tensor product $\crotimes$, where $G$ acts on $M\crotimes N$ diagonally. We can define an analogous tensor product $\crotimes$ on $\mathbf{CoSh}(\RG)$ with $G$ acting diagonally, which makes $\mathbf{CoSh}(\RG)$ symmetric monoidal. Then the adjunction between $\mathbf{PMod}(\RG)$ and $\mathbf{CoSh}(\RG)$ of Proposition \ref{propff} is still strong monoidal.
\item\label{rm3} Taking $R=\widehat{\Z}$ in \ref{rm1}, we obtain the symmetric monoidal category $(\coshZ,\widehat{\otimes}_{\widehat{\,\Z}})$. This allows us to talk about monoid objects and module objects in $\coshZ$. In fact, a profinite ring $R$ is a special instance of a monoid object in $\coshZ$, and a bundle of profinite $R$-modules is a special instance of an $R$-module object in $\coshZ$, almost by definition. Observe that if $(M,p,X)\in\coshR$ has its $R$-module structure given by a map $\rho\colon M\,\,\widehat{\otimes}_{\widehat{\,\Z}}\,\, R\to M$ in $\coshZ$, then the $R$-module structure on $\bigcoplus_X M_x$ is given by the composite $$\left(\bigcoplus_X M_x\right)\,\,\widehat{\otimes}_{\widehat{\,\Z}}\,\, R=\bigcoplus_X(M_x\,\,\widehat{\otimes}_{\widehat{\,\Z}}\,\,R)\to\bigcoplus_X M_x,$$ where the second map is obtained by applying the profinite direct sum functor to $\rho$.
\item\label{rm4} The following isomorphism given in Proposition \ref{tensorcomm} \begin{eqnarray}\label{eqn1}
\left(\bigcoplus_X M_x\right)\crotimes\left(\bigcoplus_YN_y\right)=\bigcoplus_{X\times Y}(M_x\crotimes N_y)
\end{eqnarray}
is natural in the bundles $M$ and $N$. It follows formally that if $(M,p,X)$ is a bundle of $S'$-$R$-bimodules and $(N,q,Y)$ is a bundle of $R$-$S$-bimodules, then the isomorphism in (\ref{eqn1}) is automatically an $S'$-$S$-bimodule isomorphism. Indeed, the $S'$-module structure on $M$ is given by a right $R$-module bundle map $S'\,\,\widehat{\otimes}_{\widehat{\,\Z}}\,\, M\to M$, and similarly the $S$-module structure on $N$ is given by a left $R$-module bundle map $N\,\,\widehat{\otimes}_{\widehat{\,\Z}}\,\, S\to N$. Applying the naturality of (\ref{eqn1}) to these maps gives what we claimed.
\end{enumerate}
\end{rmk}

\bibliographystyle{unsrt}
%\bibliography{bib}

\begin{thebibliography}{10}

\bibitem{melfree}
Oleg~Vladimirovich Mel'nikov.
\newblock Subgroups and the homology of free products of profinite groups.
\newblock {\em Izv. Akad. Nauk SSSR Ser. Mat.}, 53(1):97--120, 1989.

\bibitem{ribesgraph}
Luis Ribes.
\newblock {\em Profinite graphs and groups}, volume~66.
\newblock Springer, 2017.

\bibitem{gareth}
Gareth Wilkes.
\newblock Pontryagin duality and sheaves of profinite modules.
\newblock {\em J. Algebra}, 690:379--399, 2026.

\bibitem{boggi}
Marco Boggi.
\newblock Lannes' {$T$}-functor and {${\rm mod}\text{-}p$} cohomology of
  profinite groups.
\newblock {\em J. Algebra}, 694:109--147, 2026.

\bibitem{johnmodular}
J.~W. MacQuarrie.
\newblock Modular representations of profinite groups.
\newblock {\em J. Pure Appl. Algebra}, 215(5):753--763, 2011.

\bibitem{profinite}
Luis Ribes and Pavel Zalesskii.
\newblock {\em Profinite groups}.
\newblock Springer, 2000.

\bibitem{maccat}
Saunders Mac~Lane.
\newblock {\em Categories for the working mathematician}, volume~5.
\newblock Springer Science \& Business Media, 2013.

\bibitem{macsheaves}
Saunders Mac~Lane and Ieke Moerdijk.
\newblock {\em Sheaves in geometry and logic: A first introduction to topos
  theory}.
\newblock Springer Science \& Business Media, 2012.

\bibitem{weibel}
Charles~A. Weibel.
\newblock {\em An introduction to homological algebra}, volume~38 of {\em
  Cambridge Studies in Advanced Mathematics}.
\newblock Cambridge University Press, Cambridge, 1994.

\bibitem{brumer}
Armand Brumer.
\newblock Pseudocompact algebras, profinite groups and class formations.
\newblock {\em J. Algebra}, 4:442--470, 1966.

\bibitem{compactmodules}
John~W. MacQuarrie, Peter Symonds, and Pavel~A. Zalesskii.
\newblock Infinitely generated pseudocompact modules for finite groups and
  {W}eiss' theorem.
\newblock {\em Adv. Math.}, 361:106925, 35, 2020.

\bibitem{melaspherical}
Oleg~Vladimirovich Mel'nikov.
\newblock Aspherical pro-{$p$}-groups.
\newblock {\em Mat. Sb.}, 193(11):71--104, 2002.

\bibitem{garethrel}
Gareth Wilkes.
\newblock Relative cohomology theory for profinite groups.
\newblock {\em J. Pure Appl. Algebra}, 223(4):1617--1688, 2019.

\bibitem{double}
Peter Symonds.
\newblock Double coset formulas for profinite groups.
\newblock {\em Comm. Algebra}, 36(3):1059--1066, 2008.

\bibitem{talephi}
Olympia Talelli.
\newblock On groups of type {$\Phi$}.
\newblock {\em Arch. Math. (Basel)}, 89(1):24--32, 2007.

\bibitem{infinitestable}
Nadia Mazza and Peter Symonds.
\newblock The stable category and invertible modules for infinite groups.
\newblock {\em Adv. Math.}, 358:106853, 26, 2019.

\bibitem{permcom}
Peter Symonds.
\newblock Permutation complexes for profinite groups.
\newblock {\em Comment. Math. Helv.}, 82(1):1--37, 2007.

\end{thebibliography}

\end{document}